%% file: main.tex
\documentclass[conference]{IEEEtran}
\IEEEoverridecommandlockouts
\usepackage[ruled]{algorithm2e}
\usepackage{cite}
\usepackage{graphicx,epstopdf}
\usepackage{array}
\usepackage{color}
\usepackage{amsfonts}
\usepackage{amsthm,amsmath}
\graphicspath{{figures/}}
\usepackage{graphicx, pifont} 
\usepackage{tabularx}
\usepackage{lipsum}
\DeclareMathOperator{\Tr}{Tr}
\usepackage{array}
\newcolumntype{P}[1]{>{\centering\arraybackslash}p{#1}}

\begin{document}
%
\title{Health-Focused Optimal Power Flow}

\author{Logesh Kumar$^\star$, Parikshit Pareek$^\star$, Sivakumar Nadarajan, Souvik Dasgupta, Amit Gupta, and Hung D. Nguyen$^\#$
\thanks{$^\#$ Corresponding Author}
\thanks{$^\star$ The first two authors contributed equally to this work}
\thanks{Kumar is with Rolls-Royce@NTU Corporate Lab and School of EEE NTU, Singapore,  \textit{logeshkumar@ntu.edu.sg}}
\thanks{Pareek and Nguyen are with School of EEE NTU, Singapore \textit{pare0001,hunghtd@ntu.edu.sg}}
\thanks{Nadarajan, Dasgupta, and Gupta are with  Rolls-Royce Electrical, Singapore \textit{Sivakumar.Nadarajan,Souvik.Dasgupta,Amit.Gupta@Rolls-Royce.com}}
}

\maketitle


\begin{abstract}
In this paper, we propose a novel Health-Focused Optimal Power Flow (HF-OPF) to take into account the equipment health in operational and physical constraints. The health condition index is estimated based on the possible fault characteristics for generators and batteries. The paper addresses the need for understanding the relationship between health condition index and the operational constraints in OPF problems. Such a relationship is established through Finite Element Method (FEM)-based simulations. The simulations for generator and battery faults indicate the presence of an inflection point after which effect of fault becomes severe. A microgrid system is used to illustrate the performance of the proposed HF-OPF. The results show that health condition inflicts high cost of generation and can lead to infeasibility even with less critical faults.
\end{abstract}


%
\IEEEpeerreviewmaketitle

\section{Introduction}
The optimal power flow (OPF) is a fundamental problem for power system operation. It is solved for finding an economical operating point within a physical network and operational constraints \cite{low2014convex}. In recent years, the development of various microgrid architectures leads to several OPF formulations, including various sources and load types \cite{abdi2017review}. The microgrid hold conventional generators, renewable generators and as well as battery energy storage systems (BESS) \cite{levron2013optimal}. All these formulations consider all network components to be working perfectly. 
To incorporate the effect of equipment failures, an extension of classical OPF has been proposed which considered the constraints arising due to operation under a fixed list of postulated contingencies \cite{capitanescu2011state}. All these formulations handle equipment failure contingencies in binary fashion i.e., either perfectly working or completely failed. Due to this discrete handling of equipment health status, these formulations cannot be used to examine and incorporate the effects of various types of faults in generators or BESS. 
Further, Energy Management System (EMS) for on-board microgrid of an aircraft \cite{buticchi2018board}, maritime systems on ships and submarines \cite{jin2016maritime} and remote areas without grid connection \cite{mizani2009design} are required to work reliably for long and critical duration. These microgrid systems may not be able to shut down a power source as they do not have access to extra power most of the times. Also, the maintenance cannot be done on-site in most of the cases, and system may need to work properly for long-duration before it is safe to shut down for maintenance. In addition to these concerns, the advancements in online health monitoring have enabled EMS to access health conditions remotely. Hence EMS of microgrid can act upon these for preventive as well as corrective actions. Therefore, it is imperative to incorporate the effect of equipment health, in a non-discrete manner, in OPF. This will help in maintaining the system reliability with the least cost infliction. 
In this paper, we propose a novel Health-Focused Optimal Power Flow (HF-OPF) incorporating the effect of equipment health condition as a continuous function. This work attempts to understand the relationship between incipient faults and OPF constraints. A health condition index is defined for generator and BESS to reflect various fault conditions uniformly. Further, these health condition index have been converted into OPF constraints. The HF-OPF is solved using state-of-art Semi-Definite Programming (SDP) relaxation. The main contributions of this paper are:
\begin{itemize}
    \item Proposing a novel Health-Focused Optimal Power Flow (HF-OPF) formulation to incorporate equipment health by modelling OPF constraints as functions of health condition index.
    \item Defining health condition index for generator and battery for the establishment of a relationship between health and OPF constraints via FEM analysis. 
    \item Analysing the variation OPF constraints with health condition and cost inflicted by health deterioration of equipment. 
\end{itemize}

\section{Health Condition Monitoring}
For a generator, various individual faults in generator such as Turn-Turn Short Circuit (TTSC) \cite{nadarajan2014hybrid}, demagnetization fault \cite{ahsanullah2017demagnetization}, eccentricity fault have been studied independently, and individual health condition have been established. Yet, the insights into the overall, combined health condition of the generator or battery system are not well studied. This paper develops the overall health of generator and battery as a function of fault characteristics. Further, the relationship between the health condition index and OPF constraints has been established through FEM simulations. Below we elaborate the notion of health condition index in more detail.

\subsection{Health Condition Index}
For a machine, the health condition index is estimated based on the fault characteristics established in this context as fault severity, fault factor, and fault weightage. In general, calculating the health condition index of a machine is nontrivial. In this paper, we leverage the component simulations with high fidelity modelling to estimate a machine's health condition. In particular, Finite Element Method (FEM) will be performed for these fault characteristics to establish with accuracy the corresponding OPF constraints. The FEM software solves the magnetic field in case of generator to calculate the inductance from which current, voltage, and power factor are obtained. 

For a battery, the computational fluid dynamics (CFD) based ANSYS Fluent solves the electrochemistry to establish the thermal characteristics from which the electrical parameters are obtained. The obtained relationship can be used in real-time monitoring environment to detect and diagnose the health of the individual equipment.   

\begin{figure}[t]
    \centering
    \includegraphics[scale=0.25]{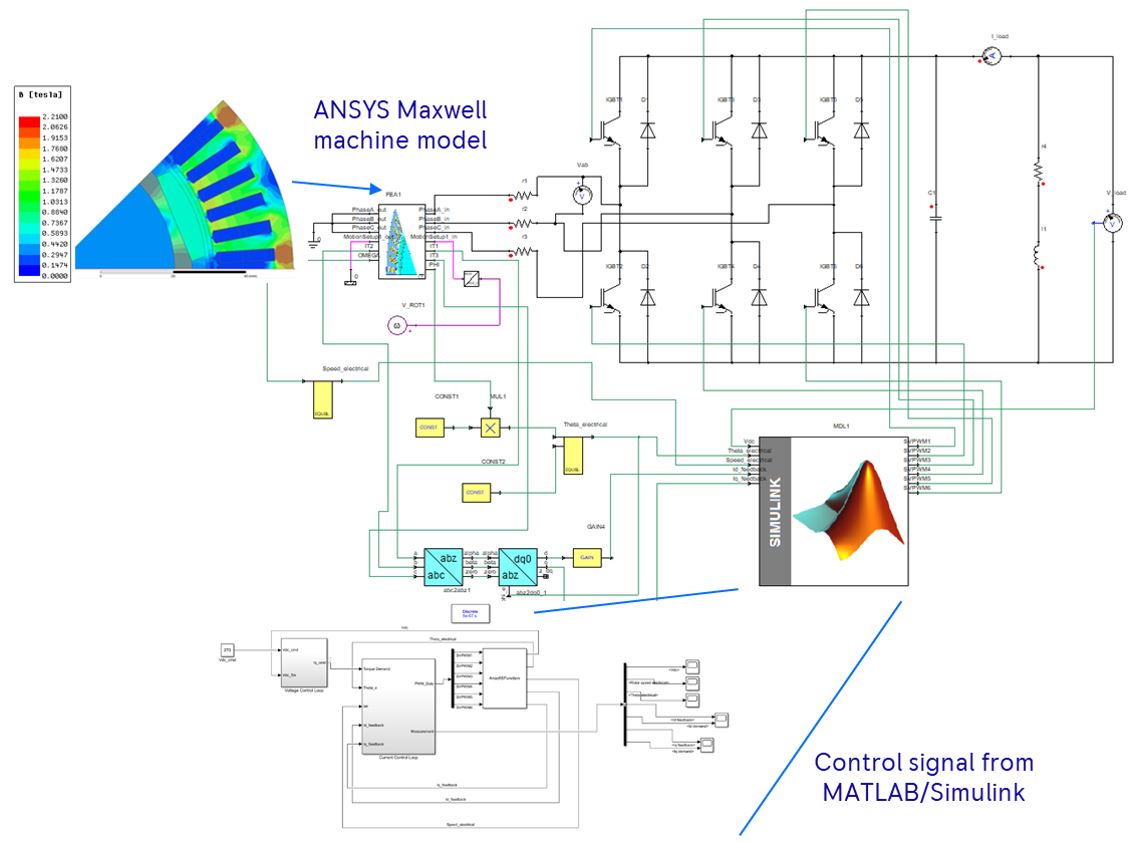}
    \vskip -0.5em
  \caption{Co-simulation of generator drive system in Twin Builder}
    \label{fig:1}
    \vskip -1.2em
\end{figure}

\begin{figure}[b]
     \centering
     \includegraphics[scale=0.44]{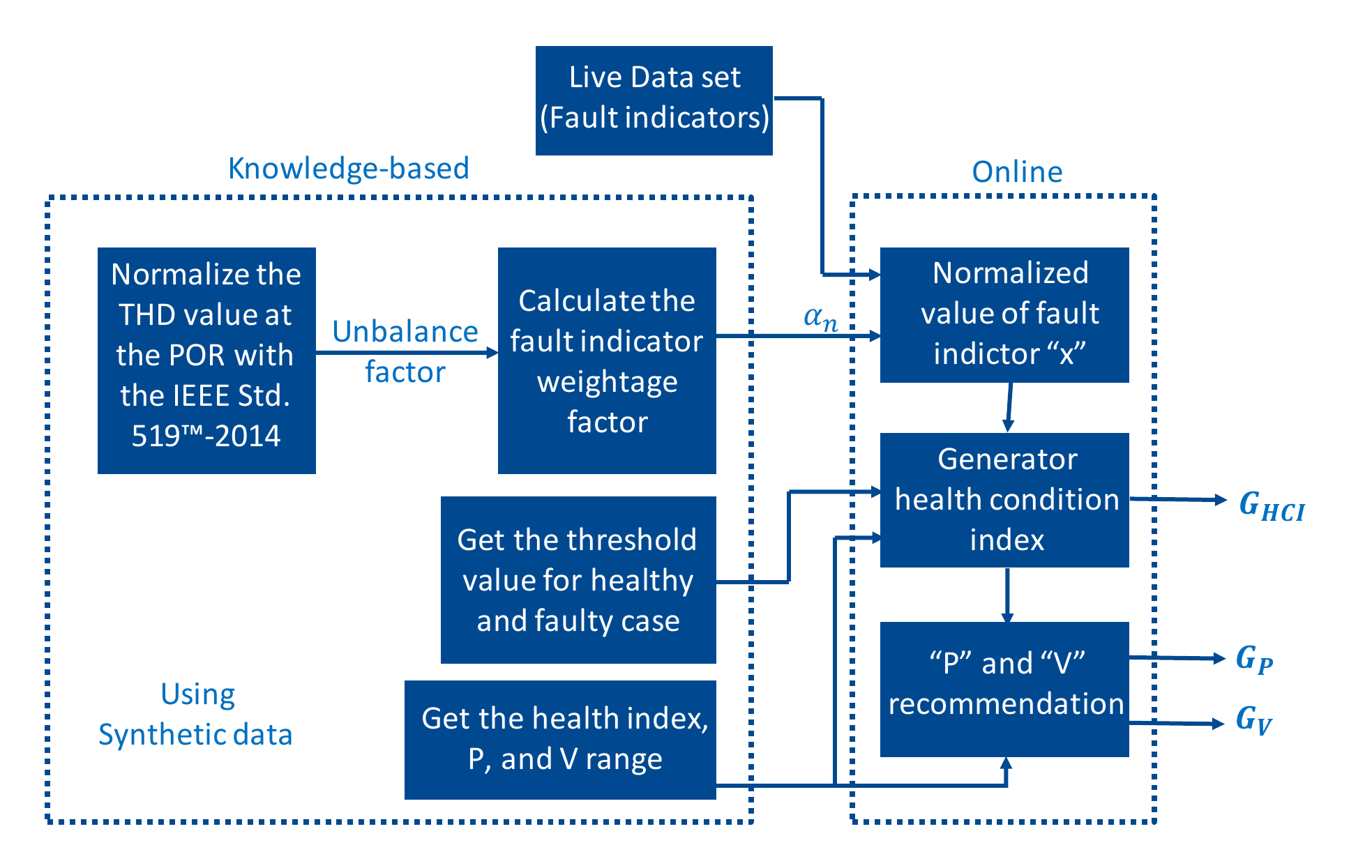}
     \vskip -0.5em
   \caption{Generator health condition index formulation block diagram}
     \label{fig:2}
     \vskip -1.2em
 \end{figure}

\subsection{Health Condition Index for Generator ($G_{HCI}$) and BESS ($B_{HCI}$)}

The generator and battery health condition index is quantified in this subsection. The key emphasis is on the definition alone, so the design aspect and their parameters are not discussed as the health condition index are normalized to their base rating. The health condition index of the generator is simulated through a generator drive system with the electromagnetic field solved in ANSYS Maxwell, and the control action is established in MATLAB-SIMULINK environment to control the fluctuating DC load voltage. The co-simulation is carried out in Twin Builder/Simplorer for the continuous exchange of data between the FEA ANSYS and MATLAB and is shown in Figure \ref{fig:1}. The electrical faults in generator include open-circuited phase fault, short-circuited fault which include TTSC, phase to phase, and phase to ground fault. The mechanical fault includes static eccentricity, dynamic eccentricity, mixed, and bearing fault. 

The TTSC faults are generally initiated by insulation breakdown caused by overheating or mechanical breakage \cite{lai2014analysis}. TTSC faults, when undetected in early-stage, may lead to other severe winding short faults due to high fault loop current. The TTSC fault produces flux in opposite direction due to the reactive nature of the fault, which increases the eddy current and hysteresis loss in the machine. The TTSC fault is simulated in ANSYS Maxwell by splitting the total conductors in coil A as a separate rectangular conductor. The separate conductors are taken externally to the Twin Builder to model it as an extra winding to carry the fault loop current. Two levels of fault severity are simulated accounting for $8.33\%$, and $25\%$. 
The demagnetization fault occurs due to either TTSC fault or prolonged overloading in the flux weakening region of the machine. In severe cases, the demagnetization and TTSC fault may lead to the saturation of the core resulting in asymmetrical flux distribution in the machine. The demagnetization fault is simulated by increasing the temperature of magnet using the BH curve. Two different fault severity of one magnet and three magnets with fault factor of $100^\circ C$ and $150^\circ C$ is simulated. 
The various causes for the mechanical faults are misalignment due to manufacturing defects, unbalanced magnetic pull, or mechanical resonance at critical speed. The eccentricity fault creates high-frequency components based on the switching frequency as well as the length of the cable. The low-frequency component is also produced in torque due to other faults in machine leading to the increase in shaft voltage and consequently to the bearing fault in the machine. In ANSYS Maxwell, the SE is introduced by shifting the stator and the coils to 72\% and 88.8\% from normal position. Similarly, DE is implemented by shifting the rotor and the magnets by 17.7\%, and 22.2\% from their initial position. 

\begin{figure}[t]
    \centering
    \vskip -1em
    \includegraphics[width=\columnwidth]{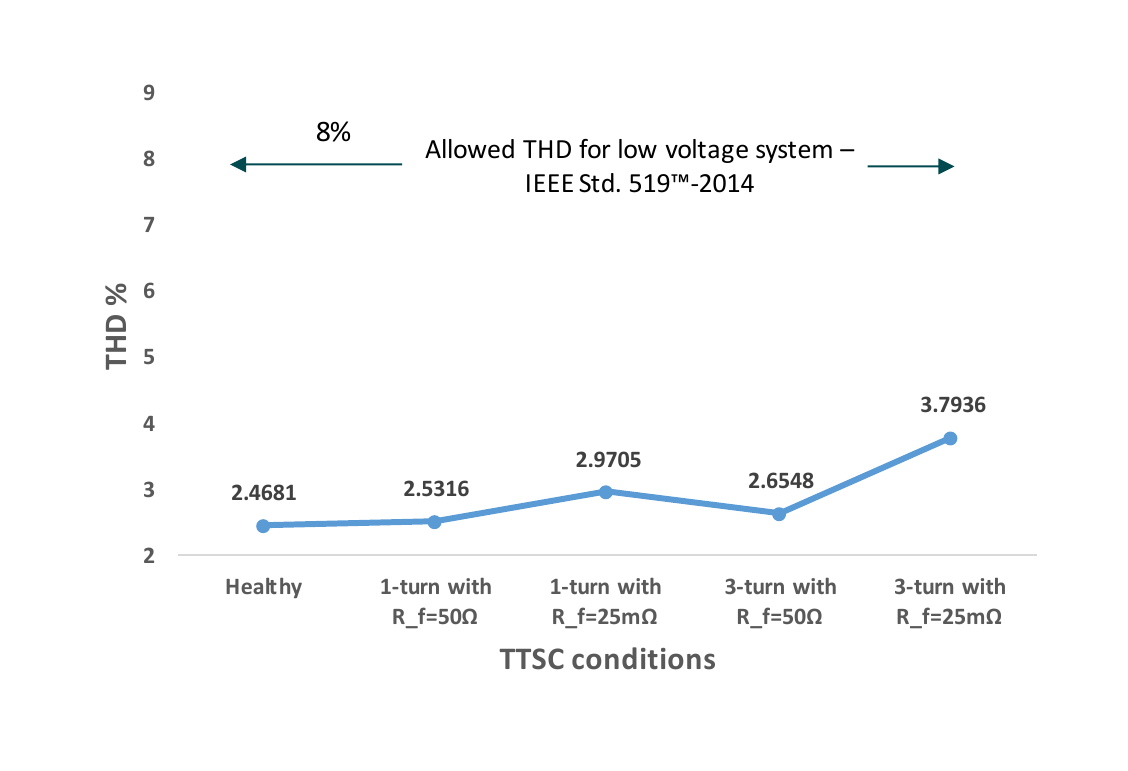}
    \vskip -2em
  \caption{Unbalanced factor calculated for TTSC fault condition}
    \label{fig:3}
     \vskip -1.2em
\end{figure}

Various faults in asset exhibit different characteristics features. Common platform to weigh these features for formulating the health index of the generator is necessary. Figure \ref{fig:2} shows the block diagram for generator health condition index, real power, and voltage formulation. The total harmonic distortion (THD) according to IEEE standard ${519^{TM}-2014}$ provides variation for any fault in the asset, so THD is used as an evaluating signature to calculate the unbalanced factor. The unbalanced factor being calculated for TTSC fault condition is shown in Figure \ref{fig:3}. Similarly, the unbalanced factor is calculated for all other generator fault based on the harmonic content. The weightage factor for the fault indicators to be used real-time are calculated based on the unbalanced factor. The real-time obtained fault indicators are normalized with their respective weightage factor. Finally, using the synthetic FEM simulation knowledge about the fault threshold, the mapping of normalized value with the health index range is done to obtain the health condition index of the generator and the values for various generator fault cases are shown in Table \ref{tab:1}. The relationship between the health condition index of the component and the constraint will not be linear. \begin{figure}[t]
    \centering
    \includegraphics[scale=0.37]{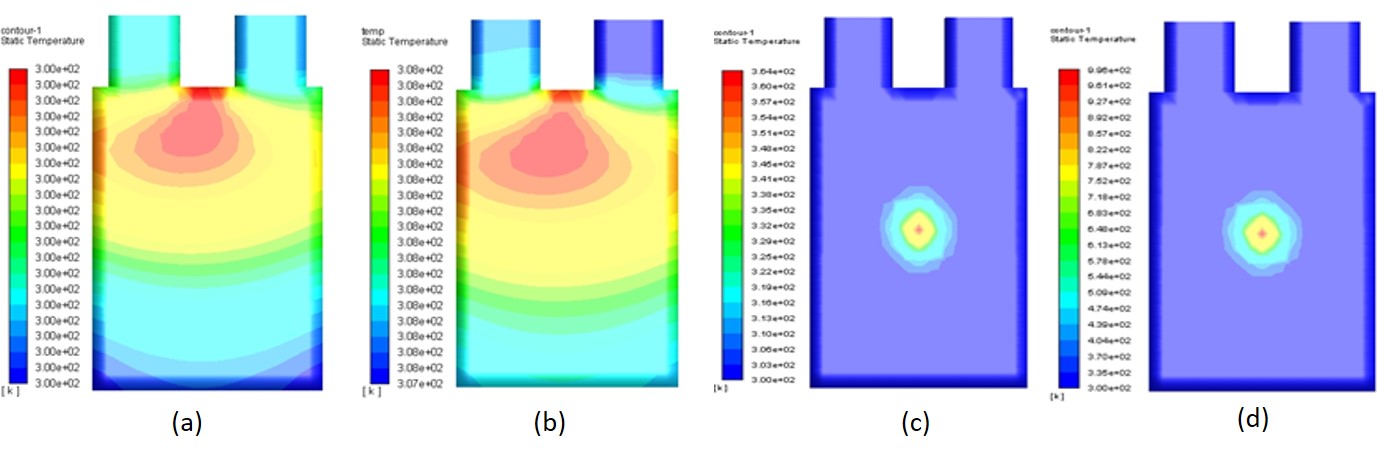}
    \vskip -0.5em
  \caption{Temperature contour with SOC of $90\%$ for (a) Healthy (b) External short-circuit (c) Internal short-circuit (d) Thermal runaway}
    \label{fig:4}
    \vskip -1.2em
\end{figure}
This is due to the definition of index and their corresponding fault characteristics. As different fault behaves differently, and their effect on system parameters vary based on the fault type and their severity. So, this non-linearity in parameters tend to occur with incipient fault cases.

\begin{table}[b]
\caption{Generator health condition index, P, and V calculated based on synthetic simulation data}
\begin{tabular}{|l|c|c|c|c|}
\hline
\multicolumn{1}{|c|}{\textbf{\begin{tabular}[c]{@{}c@{}}Generator \\ fault cases\end{tabular}}} & \textbf{\begin{tabular}[c]{@{}c@{}}Fault severity \\ based on RUL\end{tabular}} & \textbf{P range} & \textbf{V range} & \textbf{\begin{tabular}[c]{@{}c@{}}$G_{HCI}$ \\ range\end{tabular}} \\ \hline
\textbf{Healthy}                                                                                & 0                                                                               & 0                & 0                & 1                                                                \\ \hline
\textbf{Magnet fault}                                                                           & 0.5                                                                             & 0.06             & 0.075            & 0.99-0.96                                                        \\ \hline
\textbf{\begin{tabular}[c]{@{}l@{}}Static \\ eccentricity\end{tabular}}                         & 0.8                                                                             & 0.08             & 0.11             & 0.96-0.94                                                        \\ \hline
\textbf{\begin{tabular}[c]{@{}l@{}}Dynamic \\ eccentricity\end{tabular}}                        & 0.9                                                                             & 0.08             & 0.12             & 0.95-0.93                                                        \\ \hline
\textbf{\begin{tabular}[c]{@{}l@{}}Mixed \\ eccentricity\end{tabular}}                          & 0.91                                                                            & 0.07             & 0.135            & 0.92                                                             \\ \hline
\textbf{\begin{tabular}[c]{@{}l@{}}Turn-turn short-\\ circuit\end{tabular}}                     & 1                                                                               & 0.1              & 0.15             & 0.9-0.86                                                         \\ \hline
\textbf{Phase-Ground}                                                                           & 3                                                                               & 0.49             & 0.469            & 0.79                                                             \\ \hline
\textbf{Open-Phase}                                                                             & 4                                                                               & 0.59             & 0.527            & 0.76                                                             \\ \hline
\textbf{Phase-Phase}                                                                            & 5                                                                               & 0.69             & 0.587            & 0.65                                                             \\ \hline
\textbf{\begin{tabular}[c]{@{}l@{}}Three-phase \\ open\end{tabular}}                            & 7                                                                               & 0.75             & 0.827            & 0.51                                                             \\ \hline
\textbf{\begin{tabular}[c]{@{}l@{}}Bolted short-\\ circuited\end{tabular}}                      & 10                                                                              & 0.89             & 0.934            & 0                                                                \\ \hline
\end{tabular}
\label{tab:1}
\end{table}

The battery electrochemistry is solved in CFD based ANSYS Fluent. The three common faults in battery are external short circuit, internal short circuit, and thermal runaway. The internal short circuit is the incipient fault case and can lead to thermal runaway in long run. The leading causes for internal short circuit fault are lithium plating decomposition, electrode or electrolyte decomposition, or various external factors contributing to localized heat spots in the battery. The number of steps or cycles for healthy and various fault cases will be different for the same \textit{c-rate} and are represented as end-of-discharge ($F$) factor over time.
The external short circuit fault is simulated by introducing an external resistance with the battery terminal for 50\% fault severity of $0.5\Omega$. The various SOC condition for \textit{c-rate} equal to one is formulated to get $B_{HCI}$ as shown in Table \ref{tab:2}. The internal short circuit fault and thermal runaway are introduced by patching the center of the battery to a resistance of $5\times 10^{-7}\Omega$ and $5\times 10^{-8}\Omega$ respectively. The temperature contour for healthy and faulty cases with SOC of 90\% is shown in Figure \ref{fig:4}.

\begin{table}[t]
\caption{Battery health condition index, P, and V calculated based on synthetic simulation data}
\begin{tabular}{|l|c|c|c|c|}
\hline
\multicolumn{1}{|c|}{\textbf{\begin{tabular}[c]{@{}c@{}}Battery\\ fault cases\end{tabular}}} & \textbf{\begin{tabular}[c]{@{}c@{}}Fault severity \\ based on RUL\end{tabular}} & \textbf{P range} & \textbf{V range} & \textbf{\begin{tabular}[c]{@{}c@{}}$B_{HCI}$ \\ range\end{tabular}} \\ \hline
\textbf{Healthy}                                                                             & 0                                                                               & 0                & 0                & 0.9                                                              \\ \hline
\textbf{\begin{tabular}[c]{@{}l@{}}External \\ short circuit\end{tabular}}                   & 0.5                                                                             & 0.26             & 0.31             & 0.81                                                             \\ \hline
\textbf{\begin{tabular}[c]{@{}l@{}}Internal \\ short circuit\end{tabular}}                   & 0.9                                                                             & 0.91             & 0.92             & 0.53                                                             \\ \hline
\textbf{\begin{tabular}[c]{@{}l@{}}Thermal\\ runaway\end{tabular}}                           & 1                                                                               & 0.96             & 0.98             & 0                                                                \\ \hline
\end{tabular}
\label{tab:2}
\end{table}

\section{Health-Focused OPF and Its SDP Relaxation}
The conventional OPF is a constraint optimization problem for minimization of an objective while working over a set of equality and inequality constraints as:
\begin{equation}\label{eq:COPF}
    \begin{aligned}
        \min_\mathbf{x} \quad & C(\mathbf{x})\\
        \text{s.t} \quad & f_i(\mathbf{x})=0; \quad i=1,2 \dots m\\
        & x_k \leq x^u_k; \quad k=1,2 \dots n
    \end{aligned}
\end{equation}

Here, $\mathbf{x}$ is variable vector having control and state variables, $f_i(\mathbf{x})$ is equality constraint function while $x^u_k$ represents the upper limit of variable $x_k$. The upper bound vector $\mathbf{x}^u$ is considered constant if the equipment is connected. The proposed general HF-OPF formulation is given as: 
\begin{equation}\label{eq:HAOPF}
    \begin{aligned}
        \min_\mathbf{x} \quad & C(\mathbf{x})\\
        \text{s.t} \quad & f_i(\mathbf{x})=0; \quad i=1,2 \dots m\\
        & x_k \leq x^u_k(\beta_k(h_k)); \quad k=1,2 \dots n
    \end{aligned}
\end{equation}
    \vskip -0.1em
It is essential to understand that health condition index cannot be determined with certainty in various conditions. Therefore, we propose to use the probability distribution of health, $\beta_k(h_k)$, in defining the bounds on $\mathbf{x}$. This way, the HF-OPF formulation (\ref{eq:HAOPF}) will have uncertain limits which needed to take care by the chance-constrain or robust optimization tools. Currently, for better understanding of the effects of health on OPF, we use $\beta_k(h_k)=h_k$,  which can be interpreted as a mean of the probability function.  
We define network node-set as $\mathcal{N}$, generator node-set $\mathcal{G}\subseteq \mathcal{N}$ and the set of branches as $\mathcal{L}$. We define $\mathcal{B}$ as a set of nodes having BESS and use $k$ as node index. Other system parameters are defined as: 
\begin{itemize}
    \item $S_{G_k}=P_{G_k}+jQ_{G_k}$:= Apparent power generation by generator connected to $k\in \mathcal{G}$. 
    \item $S_{L_k}=P_{L_k}+jQ_{L_k}$:= Apparent power load at $k\in \mathcal{N}$.
    \item $V_k=V_{x_k}+jV_{y_k}$:= Complex node voltage at node $k\in \mathcal{N}$.
    \item $S_{kl}=P_{kl}+jQ_{kl}$:= The apparent, real, and reactive power flow in branch from node $k$ to $l$ where $(k,l) \in \mathcal{L}$.
    \item $C(P_{G_k})$:= Quadratic cost function for real power generation at node $k \in \mathcal{G}$.
    \item $\mathbf{Y}$:= Admittance matrix of size $\mathcal{N}\times\mathcal{N}$ with elements ${y}_{kl}$ where $(k,l) \in \mathcal{L}$
\end{itemize}

The ACOPF model is studied in detail by various works \cite{low2014convex,lee2019convex,nguyen2018constructing,pareek2019small}. The vital difference in proposed HF-OPF is that limits are a function of health which are considered as constant in conventional ACOPF formulation.
The HF-OPF for the real power generation cost minimization objective can be expressed as \cite{lavaei2012zero}: 
\begin{subequations}
\begin{align}
    \textrm{min} & ~~~ \sum\limits_{k\in \mathcal{G}} C(P_{G_k}) \label{eq:ACOPF_obj}\\
   \text{s.t.} \, \, & V^l_{k}(h_k) \leq |V_k| \leq V^u_{k}(h_k) \quad  \; k \in \mathcal{N}\label{eq:ACOPF_vlimit}\\
    & |S_{kl}| \leq S^u_{kl} \quad  \;(k,l) \in \mathcal{L}\label{eq:ACOPF_slimit}\\
    & P^l_{G_k}(h_k) \leq P_{G_k} \leq P^u_{G_k}(h_k) \quad \; k\in \mathcal{G}\label{eq:ACOPF_pg}\\
    & Q^l_{G_k}(h_k) \leq Q_{G_k} \leq Q^u_{G_k}(h_k) \quad \; k\in \mathcal{G}\label{eq:ACOPF_qg}\\
       &S_{G_k}=S_{D_k}+\sum \limits_{l}S_{kl}\quad \; k \in \mathcal{N}\label{eq:ACOPF_bal} \\
   & S_{kl}={y}^\star_{kl}V_kV_k^\star-{y}^\star_{kl}V_kV_l^\star \quad  \;(k,l) \in \mathcal{L} \label{eq:ACOPF_skl}
\end{align}
\end{subequations}

Here, it is essential to understand that HF-OPF considers a continuous health condition index $h_k$, which affects the capabilities of the devices connected with the network. Thus, it takes a very realistic situation into account instead of event-based analysis where only binary conditions (on/off) are considered. In the proposed formulation, the BESS is treated as a generator with lower real power limits being negative and function of health during the discharging. The additional constraint on BESS are:
\begin{subequations}
    \begin{align}
     -P^u_{k}(h_k) \leq & P_{k}(t) \leq + P^u_{k}(h_k) \quad \; k\in \mathcal{B}\label{eq:bess}\\
     E^l_{k}(h_k)\leq & E_{k}(t) \leq E^u_{k}(h_k) \quad \; k\in \mathcal{B} \label{eq:bess2}
    \end{align}
\end{subequations}

Here, $E_{k}(t)$ represents state-of-charge (SOC), in works related OPF with BESS,  $P^u_{k}$ is has been taken as rated power capacity and $E^u_{k}$/$E^l_{k}$ is a fixed number defining upper/lower limit of energy stored \cite{hossain2018integrated}. 
In the proposed HF-OPF, the complexity arises in the non-convex power flow equation (\ref{eq:ACOPF_skl}) solely due to the product term of voltage ($V_kV_l^\star$). In literature, various models have been proposed for the convex relaxation of the ACOPF problem. One of the initial model was presented in \cite{lavaei2012zero} using a lifting variable matrix $W$ which is equal to the voltage product term as $W_{kl}=V_kV^\star_l, (k,l \in \mathcal{N})$.
In this work, we use the model described in \cite{lavaei2012zero,pareek2019small}  with $W \in \mathbb{R}^{2n_b\times 2n_b}$ and $\mathbf{V}\in \mathbb{R}^{2n_b}$. The HF-OPF  problem (\ref{eq:ACOPF_obj}-\ref{eq:ACOPF_skl}), with matrix variable $W$ can be expressed as: 
\begin{subequations}
\begin{align}
    &\min  \quad \sum\limits_{k\in \mathcal{G}} C(P_{G_k}) \label{eq:WACOPF_obj}\\
   & \text{Subject to} \nonumber \\
   & P^l_{G_k}(h_k)-P_{D_k} \leq \Tr\{\mathbf{Y}_kW\} \leq P^u_{G_k}(h_k)-P_{D_k} \label{eq:WACOPF_pg}\\
     & Q^l_{G_k}(h_k)-Q_{D_k} \leq \Tr\{\mathbf{\overline{Y}}_kW\} \leq Q^u_{G_k}(h_k)-Q_{D_k} \label{eq:WACOPF_qg}\\
    & (V^l_{k}(h_k))^2 \leq \Tr \{M_kW\}\leq (V^u_{k}(h_k))^2 \label{eq:WACOPF_vlimit}\\
    & \Tr\{\mathbf{Y}_{kl}W\}^2+\Tr\{\mathbf{\overline{Y}}_{kl}W\}^2 \leq (S^u_{kl})^2 \label{eq:WACOPF_slimit}\\
    & W=\mathbf{V}\mathbf{V}^T \label{eq:WACOPF_W}
\end{align}
\end{subequations}
Here, the real and reactive power generation and $\mathbf{V}$ are:
\begin{subequations}
\begin{align}
    & P_{G_k}= \Tr\{\mathbf{Y}_kW\}+P_{D_k} \label{eq:pgrlx}\\
    & Q_{G_k}= \Tr\{\mathbf{\overline{Y}}_kW\} +Q_{D_k} \label{eq:qgrlx}\\
    & \mathbf{V}=[\mathbf{V}_x^T\;\; \mathbf{V}_y^T]^T
\end{align}
\end{subequations}
\vspace{-2pt}
The definition of the constant matrices $\mathbf{Y}_k,~\mathbf{\overline{Y}}_k,~\mathbf{Y}_{kl},~\mathbf{\overline{Y}}_{kl}$ and $M_k$ is similar the one presented in \cite{lavaei2012zero} and has not presented here to avoid repetition. 
 This non-convexity of power flow constraint (\ref{eq:ACOPF_bal}) is now transferred into a single equality constraint (\ref{eq:WACOPF_W}). For relaxed HF-OPF the equality constraint (\ref{eq:WACOPF_W}) has been replaced by: 
\begin{align}\label{eq:RCACOPF}
    W &\succeq 0.
\end{align}
The exactness of the HF-OPF has been sacrificed by removing the non-convex rank-one constraint of $W$.  The details of exactness conditions can be found in  \cite{lavaei2012zero} for similar relaxation. Thus, the SDP relaxed HF-OPF is minimization of objective \eqref{eq:WACOPF_obj} constrained to BESS SOC constrain \eqref{eq:bess2} and relaxed ACOPF constrains \eqref{eq:WACOPF_pg}- \eqref{eq:WACOPF_slimit}, \eqref{eq:pgrlx}, \eqref{eq:qgrlx} and \eqref{eq:RCACOPF} abd. The constrain \eqref{eq:WACOPF_pg} is obtained from \eqref{eq:bess} for BESS modeled as real power generator.

\section{Test Results}
\begin{figure}[t]
    \centering
    \vskip -1em
    \includegraphics[width=\columnwidth]{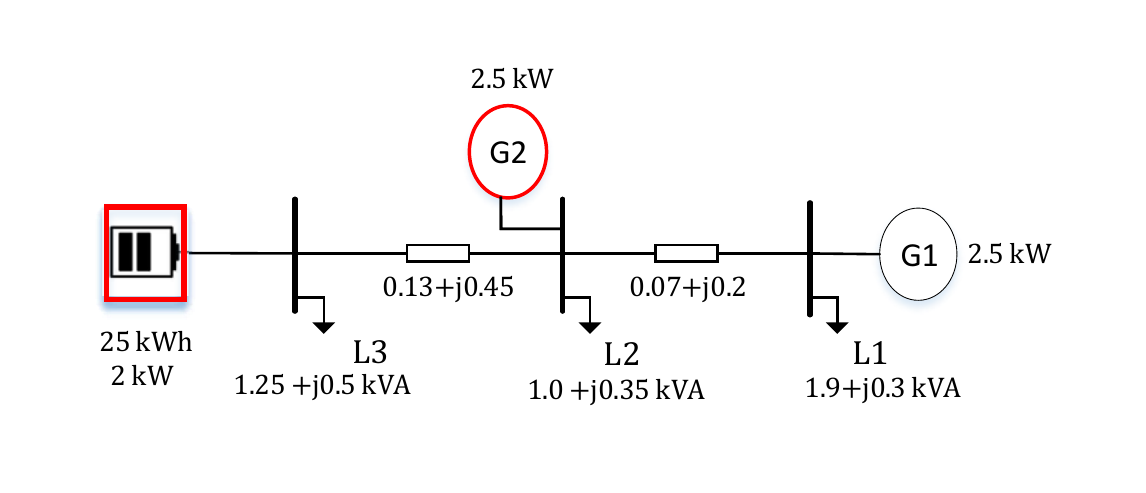}
    \vskip -1.9em
  \caption{3-Bus microgrid test network }
    \label{fig:mg}
    \vskip -1em
\end{figure}

For testing the proposed HF-OPF and analysis of the effect of health on cost, we use a 3-Bus test microgrid system as given in Fig. \ref{fig:mg}. The G1 is a costly healthy generator while different fault conditions are simulated on G2 and battery. 

First, we consider nominal loading condition. The relative cost change, with respect to ACOPF solution having with perfectly healthy generator, is calculated. Except for the variation in $P^u_2$, all other system and solver parameters are kept constant to observe only the effect of health on ACOPF solution. $P_{G_2}$ decreases sharply when the maximum output is limited due to poor health and other costly power sources (BESS and G1) output increases for load balancing. 

\begin{figure}[t]
    \centering
    \includegraphics[width=0.85\columnwidth]{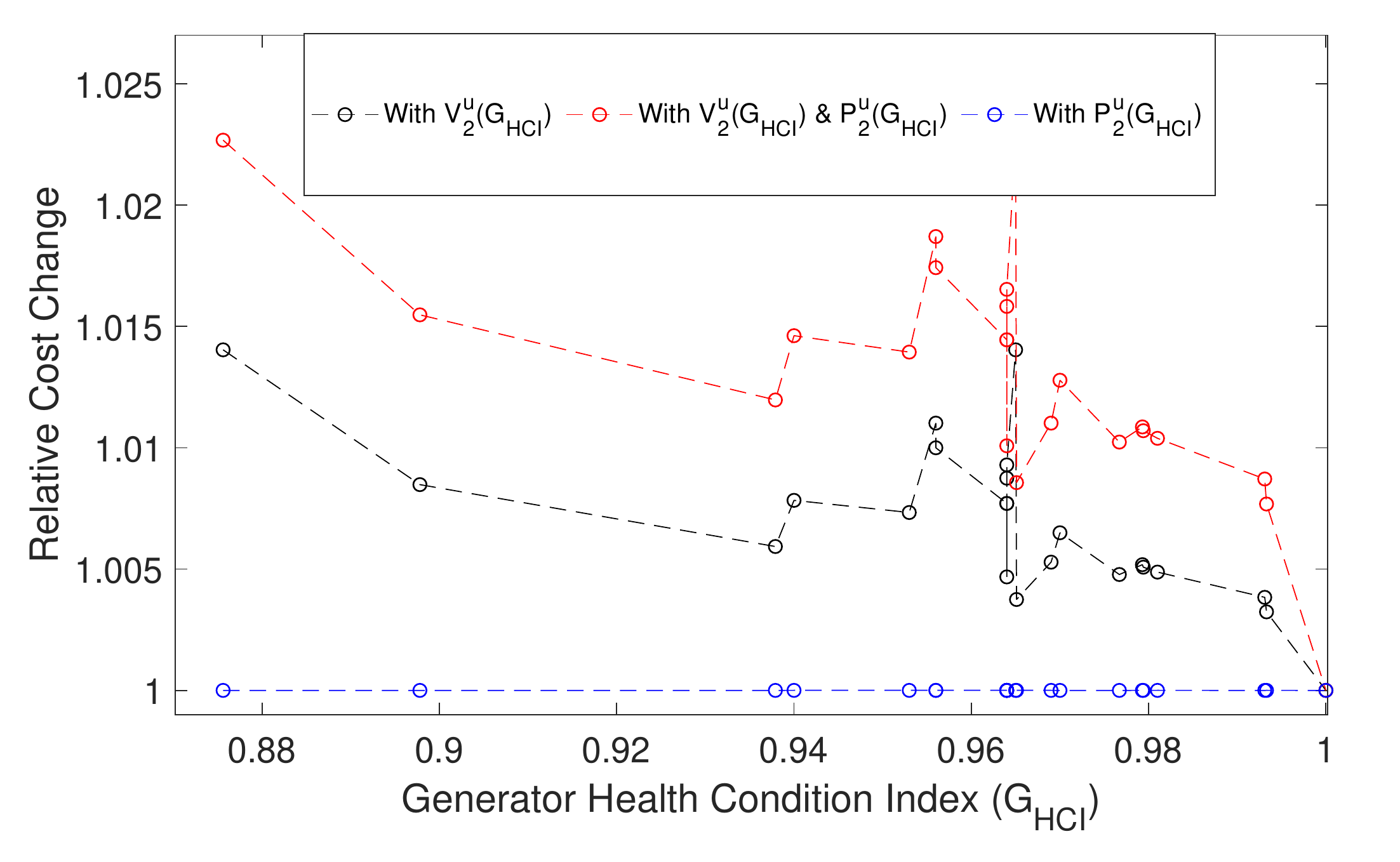}
    \vskip -1em
  \caption{Variations in relative cost change with $G_{HCI}$ at different cases}
    \label{fig:costbreak}
    \vskip -1.0em
\end{figure}

Due to the effect on voltage limit, the system may not have any solution as it will violate the voltage tolerance band (considered as $0.9-1.1$ (pu) (\ref{eq:WACOPF_vlimit})). Therefore, results of relative cost variations, including health condition's effect on voltage, are given in Figure \ref{fig:costbreak}. Considering $V^u_2$ as a function of generator health makes a more significant impact on the HF-OPF solution cost in comparison to taking only $P^u_2(G_{HCI})$ into account. These results also emphasize the requirement of HF-OPF as the effect on cost is pronounced very rapidly with a decrease of $G_{HCI}$. 

\begin{figure}[t]
    \centering
    \includegraphics[width=0.85\columnwidth]{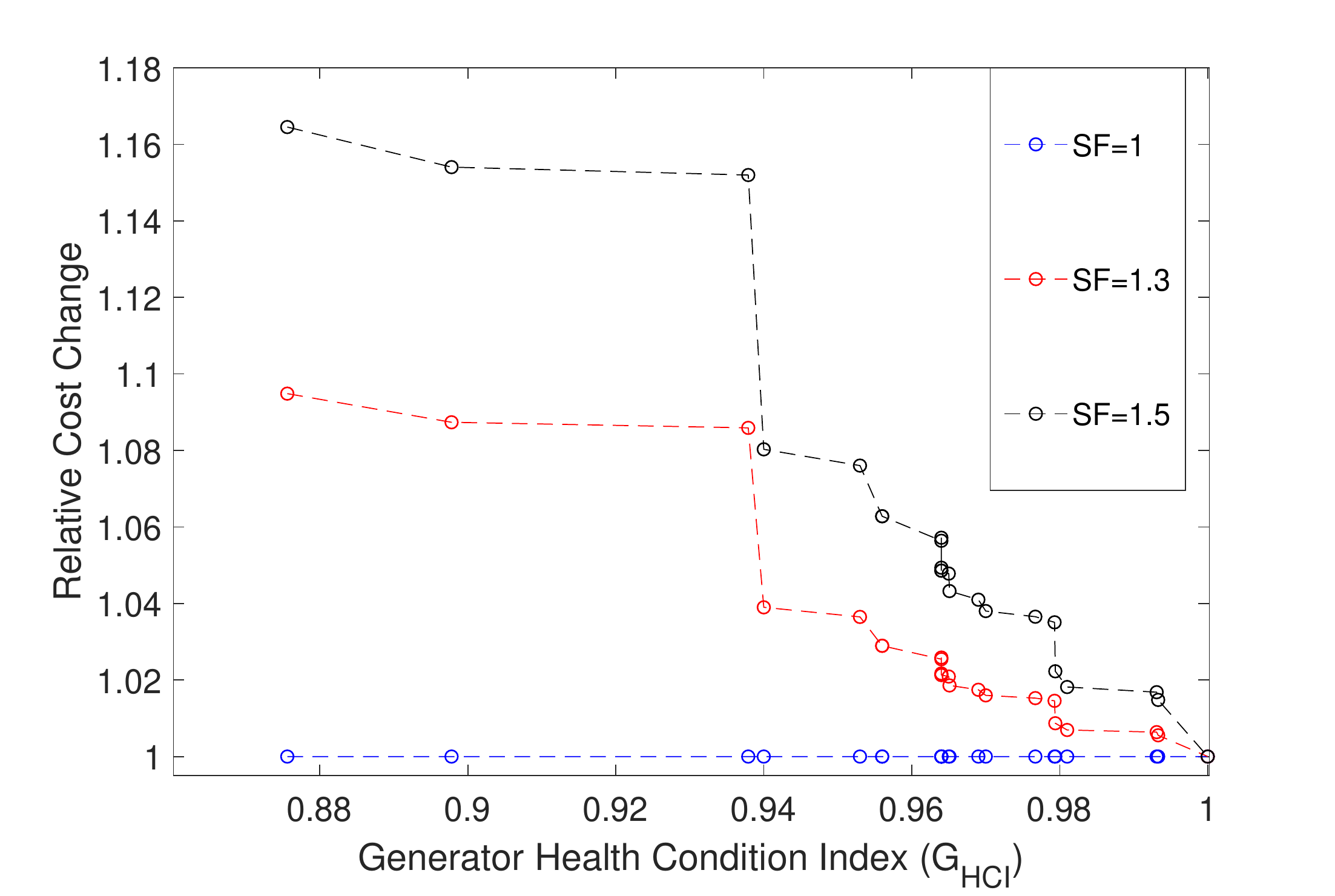}
    \vskip -1em
  \caption{Effect of loading on HF-OPF cost solution}
  \vskip -1.8em
    \label{fig:costLF}
\end{figure}

We also perform the test using three different scaling factors ($SF$). The nominal loading condition shows the same results when only $P^u_2(G_{HCI})$ is considered in Figure \ref{fig:costLF}. Subsequently, at higher loading, the effect on cost is very high with even small health decrease. These results indicate that for critical microgrid systems, the health monitoring and HF-OPF solution is essential.
Similarly, we analyze the HF-OPF results with the health condition of the battery. The cost variation is observed and step change will be similar to the generator. This variation is due to non-linearity in the effects of faults on the real power discharge capacity of the battery. These results have been able to highlight the fact that the HF-OPF solution differs from OPF solution significantly.

\section{Conclusion}
The contingency-based OPF methods consider the health of power system equipment discretely. The critical microgrid systems need to have better incorporation of equipment health to ensure high availability and reliability. To facilitate this, we propose a novel Health-Focused OPF formulation, which includes constraints as a function of health condition index. The FEM based analysis is used to establish the relationship of health condition index with constraints. The test results indicate that health issues will cause severe costs on real power generations, let alone the problems of infeasibility. The proposed HF-OPF can work with online health monitoring systems, and EMS of remote microgrid will provide improved reliability. Future works include the establishment of the nonlinear relationship between health condition index and OPF constraints as well as solving HF-OPF for uncertain health condition index and constraints. 

\section*{Acknowledgment}
This work is supported by NTU SUG and IAF-ICP, RRE Power System Health Management project.

%
\input{main.bbl}

\end{document}